\documentclass[final]{siamltex}

\usepackage{amssymb}
\usepackage{amsmath}
\usepackage{graphicx}
\usepackage{latexsym}
\usepackage{verbatim}
\usepackage{float}
\usepackage{url}
\usepackage{multirow}
\usepackage{color}
\usepackage{caption}
\usepackage{subcaption}
\usepackage{graphicx}

\newtheorem{Thm}{Theorem}

\newtheorem{Rem}{Remark}

\title{On Symmetric Positive Definite Preconditioners for Multiple Saddle-Point Systems} 

\author{John W. Pearson\thanks{School of Mathematics, The University of Edinburgh, James Clerk Maxwell Building, The King's Buildings, Peter Guthrie Tait Road, Edinburgh, EH9 3FD, United Kingdom ({\tt j.pearson@ed.ac.uk})}\and Andreas Potschka\thanks{Institute of Mathematics, Clausthal University of Technology, Erzstr. 1, 38678 Clausthal-Zellerfeld, Germany ({\tt andreas.potschka@tu-clausthal.de})}}

\begin{document}
\maketitle

\begin{abstract}
We consider symmetric positive definite preconditioners for multiple saddle-point systems of block tridiagonal form, which can be applied within the {\scshape Minres} algorithm. We describe such a preconditioner for which the preconditioned matrix has only two distinct eigenvalues, $1$ and $-1$, when the preconditioner is applied exactly. We discuss the relative merits of such an approach compared to a more widely studied block diagonal preconditioner, specify the computational work associated with applying the new preconditioner inexactly, and survey a number of theoretical results for the block diagonal case. Numerical results validate our theoretical findings.
\end{abstract}



\pagestyle{myheadings}
\thispagestyle{plain}
\markboth{J. W. PEARSON AND A. POTSCHKA}{PRECONDITIONERS FOR MULTIPLE SADDLE-POINT SYSTEMS}

\section{Introduction}\label{sec:1}

Multiple saddle-point systems and preconditioning strategies for their efficient numerical solution have attracted wide interest of late (see \cite{AliBeikBenzi,BSZ,MNN,SognZulehner}, for instance). In this paper we consider multiple saddle-point systems of the following block tridiagonal form:
\begin{equation}\label{Ak}
\mathcal{A}_k=\left(\begin{array}{cccccc}
A_0 & B_1^T & O & \hdots & O & O \\ B_1 & -A_1 & B_2^T & \ddots & & O \\ O & B_2 & A_2 & B_3^T & \ddots & \vdots \\ \vdots & \ddots & \ddots & \ddots & \ddots & O \\ O & & \ddots & B_{k-1} & (-1)^{k-1} A_{k-1} & B_k^T \\ O & O & \hdots & O & B_k & (-1)^k A_k \\
\end{array}\right),
\end{equation}
with $A_0$ symmetric positive definite, and $A_1, A_2, ..., A_k$ symmetric positive semi-definite. These conditions are not necessary to guarantee invertibility of $\mathcal{A}_k$, but are borne out of the structures of the preconditioners we wish to consider. In particular, one popular choice of preconditioner for $\mathcal{A}_k$ is the block diagonal matrix \cite{MNN,PearsonPotschka,SognZulehner}
\begin{equation}\label{BD}
\mathcal{P}_{D,k}=\left(\begin{array}{ccccc}
A_0 & O & \hdots & O & O \\ O & S_1 & \ddots & & O \\ \vdots & \ddots & S_2 & \ddots & \vdots \\ O & & \ddots & \ddots & O \\ O & O & \hdots & O & S_k
\end{array}\right),
\end{equation}
where $S_1=A_1+B_1 A_0^{-1}B_1^T$, $S_j=A_j+B_j S_{j-1}^{-1}B_j^T$ for $j=2,...,k$. Note the structure of the Schur complements $S_1,S_2,..., S_k$ indicates they are likely to be progressively more difficult to invert or approximate as $j$ increases.

{\scshape Main Assumption.} \, 
We assume for the forthcoming analysis that $\text{ker}(A_j)\cap\text{ker}(B_j^T)=\{0\}$ for $j=1,...,k$, where $\text{ker}(\cdot)$ denotes the kernel of a matrix. An important consequence of this assumption is that the Schur complements $S_1, S_2, ..., S_k$ are symmetric positive definite.

The preconditioner $\mathcal{P}_{D,k}$ has the attractive property that it may be applied within the {\scshape Minres} algorithm \cite{minres} when solving linear systems involving $\mathcal{A}_k$, with a guaranteed convergence rate based on the eigenvalues of the preconditioned matrix $\mathcal{P}_{D,k}^{-1}\mathcal{A}_k$. This can be a valuable descriptor: for instance in \cite[Theorem 4]{BrGr21} and \cite[Theorem 5.3]{PearsonPotschka} it is shown that for $k=2$, the eigenvalues of $\mathcal{P}_{D,k}^{-1}\mathcal{A}_k$ are contained within $[-\frac{1}{2}(1+\sqrt{5}),\frac{1}{2}(1-\sqrt{5})]\cup[2\,\text{cos}(\frac{3\pi}{7}),2\,\text{cos}(\frac{\pi}{7})]$, approximately $[-1.618,-0.618]\cup[0.445,1.802]$. However, the distribution of the eigenvalues is not robust with respect to $k$, as shown by the comprehensive analysis of $\mathcal{P}_{D,k}^{-1}\mathcal{A}_k$ for $k\geq1$ in the paper \cite{SognZulehner}, where the matrices are analysed as continuous operators. Indeed under certain assumptions, specifically that $A_1, A_2, ..., A_k$ are zero matrices, it is shown \cite[Theorem 2.6]{SognZulehner} that the eigenvalues of $\mathcal{P}_{D,k}^{-1}\mathcal{A}_k$ are given by $2\,\text{cos}(\frac{2i+1}{2j+3}\,\pi)$, for $j=0,1,...,k$, $i=0,1,...,j$. The eigenvalues of smallest magnitude therefore tend towards zero as $k$ increases. Hence, robust convergence of {\scshape Minres} with respect to $k$ cannot be expected when $\mathcal{P}_{D,k}$ is applied.


In this paper we propose a symmetric positive definite preconditioner for which the preconditioned matrix has only two distinct eigenvalues when it is applied exactly, independently of the number of blocks $k$, and which in general involves approximating the same matrices $A_0, S_1, S_2, ..., S_k$ as does $\mathcal{P}_{D,k}$. 
We showcase the benefit of the new preconditioner with practical computations, where the number of approximate inversions of the final (and typically computationally most expensive) Schur complement can frequently be cut in half.
We rigorously analyse the eigenvalue clustering for the case where the final Schur complement $S_k$ is approximated. We do not present a result for the much more difficult case in which all Schur complements are approximated, which is a highly challenging task even for the simpler preconditioner $\mathcal{P}_{D,k}$ in the case $k=2$ (see \cite{BrGr21}).
While we focus on the block tridiagonal structure \eqref{Ak}, it is likely that this approach could be generalized to multiple saddle-point systems with larger bandwidth. 

After introducing the new preconditioner, providing theoretical results, and discussing implementation details in Section \ref{sec:2}, we compare its numerical properties with those of the preconditioner $\mathcal{P}_{D,k}$ in Section \ref{sec:3}, and present results of tests on two applications from PDE-constrained optimization in Section \ref{sec:4}.

{\scshape Notation.}
Our definition of $\mathcal{A}_k$ involves $k+1$ matrix blocks: $k=1$ therefore corresponds to classical (generalized) saddle-point systems, $k=2$ to double saddle-point systems, $k=3$ to triple saddle-point systems, and so on. Let the $j$-th diagonal block of $\mathcal{A}_k$ contain $n_{j-1}$ rows and columns, so $\mathcal{A}_k$ has dimension $\sum_{j=0}^k n_j$. We denote by $I_n$ an identity matrix of dimension $n$, and by $O$ a zero matrix of appropriate size. As above, $\mathcal{P}_{D,k}$ is a block diagonal preconditioner of the form \eqref{BD} defined for any integer $k\geq1$, and $\mathcal{P}_{L,k}$ and $\mathcal{P}_{U,k}$ are analogous block lower triangular and block upper triangular matrices to be defined.

\section{An Ideal Symmetric Positive Definite Preconditioner}\label{sec:2}

We now suggest a modified preconditioner, defined for any integer $k\geq1$:
\begin{equation*}
\mathcal{P}_k=\mathcal{P}_{L,k}\mathcal{P}_{D,k}^{-1}\mathcal{P}_{U,k},
\end{equation*}
where
\begin{equation*}
\mathcal{P}_{L,k}=\left(\begin{array}{ccccc}
A_0 & O & \hdots & O & O \\ B_1 & -S_1 & \ddots & & O \\ O & B_2 & S_2 & \ddots & \vdots \\ \vdots & \ddots & \ddots & \ddots & O \\ O & \hdots & O & B_k & (-1)^k S_k \\
\end{array}\right)
\end{equation*}
and $\mathcal{P}_{U,k}=\mathcal{P}_{L,k}^T$. Note $\mathcal{P}_k$ is clearly symmetric positive definite by the Main Assumption of Section \ref{sec:1}. Furthermore the eigenvalues of $\mathcal{P}_k^{-1}\mathcal{A}_k$ have a very convenient structure, as we prove below.

\subsection{Theoretical results}

The following is a general result about the eigenvalues of the preconditioned matrix $\mathcal{P}_k^{-1}\mathcal{A}_k$, for our new preconditioner:
\begin{Thm}\label{thm:Multiple}
For any multiple saddle-point system ($k\geq1$), the preconditioned matrix $\mathcal{P}_k^{-1}\mathcal{A}_k$ has two distinct eigenvalues. Specifically, it has $\sum_{j=0}^{\lfloor k/2 \rfloor}n_{2j}$ eigenvalues equal to $1$, and $\sum_{j=1}^{\lfloor (k+1)/2 \rfloor}n_{2j-1}$ eigenvalues equal to $-1$.
\proof Straightforward linear algebra tells us that
\begin{align}\label{PLA}
\mathcal{P}_{L,k}^{-1}\mathcal{A}_k=\left(\begin{array}{cccccc}
I_{n_0} & S_0^{-1}B_1^T & O & \hdots & O & O \\ O & I_{n_1} & -S_1^{-1}B_2^T & \ddots & & O \\ \vdots & \ddots & I_{n_2} & S_2^{-1}B_3^T & \ddots & \vdots \\ \vdots & & \ddots & \ddots & \ddots & O \\ O & & & O & I_{n_{k-1}} & (-1)^{k-1}S_{k-1}^{-1}B_k^T \\ O & O & \hdots & \hdots & O & I_{n_k} \\
\end{array}\right), \\
\label{PDPU} \mathcal{P}_{D,k}^{-1}\mathcal{P}_{U,k}=\left(\begin{array}{cccccc}
I_{n_0} & S_0^{-1}B_1^T & O & \hdots & O & O \\ O & -I_{n_1} & S_1^{-1}B_2^T & \ddots & & O \\ \vdots & \ddots & I_{n_2} & S_2^{-1}B_3^T & \ddots & \vdots \\ \vdots & & \ddots & \ddots & \ddots & O \\ O & & & O & (-1)^{k-1}I_{n_{k-1}} & S_{k-1}^{-1}B_k^T \\ O & O & \hdots & \hdots & O & (-1)^k I_{n_k} \\
\end{array}\right),
\end{align}
where $S_0:=A_0$. Using \eqref{PDPU}, it holds that $\mathcal{P}_{U,k}^{-1}\mathcal{P}_{D,k}=(\mathcal{P}_{D,k}^{-1}\mathcal{P}_{U,k})^{-1}$ is an upper triangular matrix with block diagonal entries alternating between identity matrices (of dimensions $n_0,n_2,...$) and negative identity matrices (of dimensions $n_1,n_3,...$). Using this information, along with \eqref{PLA}, tells us that $\mathcal{P}_k^{-1}\mathcal{A}_k=(\mathcal{P}_{U,k}^{-1}\mathcal{P}_{D,k})(\mathcal{P}_{L,k}^{-1}\mathcal{A}_k)$ is also block upper triangular, with the same block diagonal entries as $\mathcal{P}_{U,k}^{-1}\mathcal{P}_{D,k}$. For example, for $k=2$,
\begin{equation*}
\mathcal{P}_2^{-1}\mathcal{A}_2=(\mathcal{P}_{U,2}^{-1}\mathcal{P}_{D,2})(\mathcal{P}_{L,2}^{-1}\mathcal{A}_2)=\left(\begin{array}{ccc}
I_{n_0} & 2S_0^{-1}B_1^T & -2S_0^{-1}B_1^T S_1^{-1}B_2^T \\ O & -I_{n_1} & 2S_1^{-1}B_2^T \\ O & O & I_{n_2} \\
\end{array}\right).
\end{equation*}
Reading off block diagonal entries of the triangular matrix $\mathcal{P}_k^{-1}\mathcal{A}_k$ gives the result. $\square$
\end{Thm}

Theorem \ref{thm:Multiple} tells us that applying {\scshape Minres} to solve a linear system involving the matrix $\mathcal{A}_k$ (for any $k\geq1$), with preconditioner $\mathcal{P}_k$, leads to guaranteed convergence in at most two iterations in exact arithmetic.

As highlighted above, the final Schur complement $S_k$ may be expected to frequently be the Schur complement of most complex structure, and hence the most difficult to approximate. It is therefore illustrative to highlight a result when a suitable approximation to $S_k$ is applied within the preconditioner.

\begin{Thm}\label{thm:ApproxSk}
Let $\widehat{\mathcal{P}}_k$ be the preconditioner $\mathcal{P}_k$, but with $S_k$ replaced by a symmetric positive definite approximation $\widehat{S}_k$. Then for any multiple saddle-point system ($k\geq1$), the eigenvalues of the preconditioned matrix $\widehat{\mathcal{P}}_k^{-1}\mathcal{A}_k$ are equal to $1$, $-1$, or are within the range of eigenvalues of $(-1)^k\widehat{S}_k^{-1}S_k$.
\proof This result may be obtained by applying identical analysis to Theorem \ref{thm:Multiple}: $\mathcal{P}_{U,k}^{-1}\mathcal{P}_{D,k}$ and $\mathcal{P}_{L,k}^{-1}\mathcal{A}_k$ have the same structure, except the bottom-right entry of $\mathcal{P}_{L,k}^{-1}\mathcal{A}_k$ is given by $\widehat{S}_k^{-1}S_k$. Reading off block diagonal entries of $\widehat{\mathcal{P}}_k^{-1}\mathcal{A}_k$ gives the result, specifically that it has $\sum_{j=0}^{\lfloor (k-1)/2 \rfloor}n_{2j}$ eigenvalues equal to $1$, $\sum_{j=1}^{\lfloor k/2 \rfloor}n_{2j-1}$ eigenvalues equal to $-1$, with the remaining $n_k$ eigenvalues given by those of $(-1)^k\widehat{S}_k^{-1}S_k$. $\square$
\end{Thm}

Theorem \ref{thm:ApproxSk} highlights that the eigenvalues of the preconditioned matrix remain tightly clustered when the final Schur complement $S_k$ is well approximated. Of course, the effective approximation of $S_k$ is an important mathematical question, in particular to ensure that the eigenvalues of $\widehat{S}_k^{-1}S_k$ (and hence those of the preconditioned system) do not become close to zero or very large in magnitude. The design of such an approximation is a problem-specific question, and should be answered to ensure either $\widehat{\mathcal{P}}_{D,k}$ or $\widehat{\mathcal{P}}_k$ is effective in practice. See \cite[Theorem 6]{BrGr21}, which indicates how the approximation of $S_k$, as well as that of the previous blocks, impacts the eigenvalues of the preconditioned system with the block diagonal preconditioner in the $k=2$ case. This analysis is complicated even for smaller values of $k$, suggesting that the corresponding analysis of the preconditioner $\widehat{\mathcal{P}}_k$ with inexact approximations of the blocks would be an interesting and challenging research direction. We highlight that, for many problems, it is possible to design mesh-robust approximations of all Schur complements, in particular $S_k$, mitigating the issue of outlier eigenvalues and leading to rapid convergence when using $\widehat{\mathcal{P}}_{D,k}$ or $\widehat{\mathcal{P}}_k$: see \cite[Section 5.2]{BrGr21}, \cite[Section 8]{MNN}, \cite[Section 3]{SognZulehner}, and Section \ref{sec:4-1} of this paper, for examples of such problems.

We now make some further observations about the above results.


\begin{Rem}
Applying Theorem \ref{thm:Multiple} with $k=1$ tells us that when preconditioning a classical (generalized) saddle-point system with
\begin{equation*}
\mathcal{P}_1=\left(\begin{array}{cc}
A_0 & O \\ B_1 & -S_1 \\
\end{array}\right)\left(\begin{array}{cc}
I_{n_0} & A_0^{-1}B_1^T \\ O & -I_{n_1} \\
\end{array}\right)=\left(\begin{array}{cc}
A_0 & B_1^T \\ B_1 & S_1+B_1 A_0^{-1}B_1^T \\
\end{array}\right),
\end{equation*}
the preconditioned system $\mathcal{P}_1^{-1}\mathcal{A}_1$ has $n_0$ eigenvalues equal to $1$, and $n_1$ eigenvalues equal to $-1$. This is highlighted in \cite{GMPS,MGW} for $A_1=O$. Note that in this case $\mathcal{P}_{D,k}$ itself ensures that the preconditioned matrix has only three distinct eigenvalues, namely $1$ and $\frac{1}{2}(1\pm\sqrt{5})$ \cite{Kuznetsov,MGW}. As $k$ increases, however, the clustering properties with the block diagonal preconditioner become less favorable, while $\mathcal{P}_k$ is robust with respect to $k$ when applied exactly.
\end{Rem}


\begin{Rem}
The matrix $\mathcal{A}_k$ can be decomposed as $\mathcal{A}_k=\mathcal{P}_{L,k}\bar{\mathcal{P}}_{D,k}^{-1}\mathcal{P}_{U,k}$, where
\begin{equation*}
\bar{\mathcal{P}}_{D,k}=\left(\begin{array}{ccccc}
A_0 & O & \hdots & O & O \\ O & -S_1 & \ddots & & O \\ \vdots & \ddots & S_2 & \ddots & \vdots \\ O & & \ddots & \ddots & O \\ O & O & \hdots & O & (-1)^k S_k
\end{array}\right).
\end{equation*}
The structures of $\mathcal{A}_k$ and $\mathcal{P}_k=\mathcal{P}_{L,k}\mathcal{P}_{D,k}^{-1}\mathcal{P}_{U,k}$, as well as the eigenvalue result in Theorem \ref{thm:Multiple}, inform us that although $\mathcal{P}_k$ should be an excellent preconditioner for $\mathcal{A}_k$, it is not close to an actual block decomposition of the original matrix, even in the ideal case. Indeed, $\mathcal{A}_k$ is an indefinite matrix, whereas $\mathcal{P}_k$ is positive definite by design, so that it may be accommodated within {\scshape Minres}. Any preconditioner which results in a very close representation of $\mathcal{A}_k^{-1}$ would in general need to be applied within a non-symmetric iterative method such as {\scshape Gmres} \cite{gmres}.
\end{Rem}

\begin{Rem}
Equation \eqref{PLA} informs us that using $\mathcal{P}_{L,k}$ (or $\mathcal{P}_{U,k}$) as a preconditioner leads to a preconditioned matrix with only one distinct eigenvalue, specifically $1$. A suitable non-symmetric iterative method such as {\scshape Gmres} could therefore be applied. 
However, the preconditioned matrix in this case would be derogatory: even in the most convenient setting of equally-sized blocks, the minimum polynomial would in general be $(x-1)^{k+1}$, so convergence would not be independent of $k$ if the preconditioner were exactly applied (see \cite[Proposition 2 (4)]{gmres}). 
Furthermore, convergence of non-symmetric iterative methods cannot in general be described using eigenvalues, so using $\mathcal{P}_{L,k}$ (or $\mathcal{P}_{U,k}$) as a preconditioner with approximations of the diagonal blocks would not give us a guaranteed convergence rate, unlike {\scshape Minres} with $\mathcal{P}_k$ or $\mathcal{P}_{D,k}$. 
In \cite[Section 3]{CJL22}, as well as considering $\mathcal{P}_{L,k}$ the authors examine preconditioning $\mathcal{A}_k$ with $\bar{\mathcal{P}}_{D,k}$. In our setting, with the additional assumption that each $A_j$, $j=0,1,...,k$, is symmetric positive definite, the authors show that the eigenvalues of the preconditioned system are positive real; the authors also analyse situations where $\mathcal{A}_k$ is non-symmetric, if $A_j=O$, $j=1,...,k$. In constrast to the approach of this paper, however, the fact that a preconditioner of the form $\bar{\mathcal{P}}_{D,k}$ is not positive definite means that it may also not be applied within {\scshape Minres}.
\end{Rem}

\subsection{Implementation}

To ensure a practical solver, one would not apply the inverses of the matrices $S_0, S_1, ..., S_k$ exactly to relevant vectors. One would instead replace these with suitable symmetric positive definite approximations $\widehat{S}_0, \widehat{S}_1, ..., \widehat{S}_{k-1}$, as well as $\widehat{S}_k$ as in Theorem \ref{thm:ApproxSk}, within $\mathcal{P}_{L,k}$, $\mathcal{P}_{D,k}$, $\mathcal{P}_{U,k}$, where $\widehat{S}_0:=\widehat{A}_0$. Let us again denote the corresponding preconditioner $\widehat{\mathcal{P}}_k=\widehat{\mathcal{P}}_{L,k}\widehat{\mathcal{P}}_{D,k}^{-1}\widehat{\mathcal{P}}_{U,k}$, which remains symmetric positive definite. The most efficient way of applying this preconditioner is to apply the inverse of
\begin{equation*}
\underbrace{\left(\begin{array}{ccccc}
\widehat{S}_0 & O & \hdots & O & O \\ B_1 & -\widehat{S}_1 & \ddots & & O \\ O & B_2 & \widehat{S}_2 & \ddots & \vdots \\ \vdots & \ddots & \ddots & \ddots & O \\ O & \hdots & O & B_k & (-1)^k \widehat{S}_k \\
\end{array}\right)}_{\widehat{\mathcal{P}}_{L,k}}\underbrace{\left(\begin{array}{ccccc}
I_{n_0} & \widehat{S}_0^{-1}B_1^T & O & \hdots & O \\ O & -I_{n_1} & \widehat{S}_1^{-1}B_2^T & \ddots & \vdots \\ \vdots & \ddots & I_{n_2} & \ddots & O \\ O & & \ddots & \ddots & \widehat{S}_{k-1}^{-1}B_k^T \\ O & O & \hdots & O & (-1)^k I_{n_k} \\
\end{array}\right)}_{\widehat{\mathcal{P}}_{D,k}^{-1}\widehat{\mathcal{P}}_{U,k}}.
\end{equation*}
Given the above formulation for applying (the inverse of) $\widehat{\mathcal{P}}_k$ in practice, we make some observations about the computational cost of the preconditioner, and discuss when it is likely to be most effective compared to $\widehat{\mathcal{P}}_{D,k}$.

\textbf{Computational cost of applying $\widehat{\mathcal{P}}_k^{-1}$.} \, Examining the above decomposition of $\widehat{\mathcal{P}}_k$, we observe that the computational cost of applying the inverse preconditioner involves approximating the inverse operations of $\widehat{S}_0, \widehat{S}_1, ..., \widehat{S}_{k-1}$ twice per {\scshape Minres} iteration, and the inverse of the final Schur complement approximation $\widehat{S}_k$ only once. Note that the block diagonal preconditioner $\widehat{\mathcal{P}}_{D,k}$ requires only one application of each block per iteration. However, as shown in Theorems \ref{thm:Multiple} and \ref{thm:ApproxSk}, a key advantage of this approach is that the symmetric positive definite preconditioner $\widehat{\mathcal{P}}_k$ can lead to improved clustering of eigenvalues compared to $\widehat{\mathcal{P}}_{D,k}$, while requiring no additional applications of the most complicated Schur complement.

\textbf{Parallelism.} \, One advantage of the preconditioner $\mathcal{P}_{D,k}$ as opposed to $\mathcal{P}_k$ studied here is that the inverse operators of the individual blocks of the preconditioner $S_0, S_1, ..., S_k$ can be approximated in parallel over the number of blocks $k$. There is no obvious way of doing this for $\mathcal{P}_k$, due to the block triangular matrices involved, although of course any available parallelism may of course be exploited \emph{within} the blocks. If parallelization is not a focus of a solver, however, the preconditioner $\mathcal{P}_k$ is likely to perform better on a multiple saddle-point system than $\mathcal{P}_{D,k}$, due to the improved clustering of the eigenvalues of the preconditioned matrix.

\textbf{Potential effectiveness of $\mathcal{P}_k$.} \, The preconditioner $\mathcal{P}_k$ has two valuable properties when applied to multiple saddle-point systems: it is symmetric positive definite, which allows provable convergence bounds for {\scshape Minres} based on the preconditioned matrix, and has a small number (only two) of distinct eigenvalues of this preconditioned matrix when the `ideal' preconditioner is applied. It uses the same approximation steps as $\mathcal{P}_{D,k}$, specifically approximating the inverse operations of $S_0, S_1, ..., S_k$, with two applications of all but the final Schur complement required per {\scshape Minres} iteration. The approach is particularly effective if $S_0, S_1, ..., S_{k-1}$ may be well approximated in a relatively low CPU time, in particular compared to $S_k$. This is partly due to the fact that $S_0, S_1, ..., S_{k-1}$ must be approximated twice per iteration within $\mathcal{P}_k$ (as opposed to once within $\mathcal{P}_{D,k}$), and partly because earlier Schur complements (while less complex in structure) may have a detrimental impact if poorly approximated, due to their propagation through the system. The final Schur complement $S_k$ is frequently very complex but may have less influence on the overall quality of the preconditioner. Hence, inaccurate approximations of $S_0, S_1, ..., S_{k-1}$ could be amplified to a greater extent within $\mathcal{P}_k$ than $\mathcal{P}_{D,k}$, however given their frequently simpler structure it is worthwhile to invest computational effort in approximating these terms. The gain in performance when applying $\mathcal{P}_k$ instead of $\mathcal{P}_{D,k}$ is hence likely to be maximized for larger $k$, when $S_k$ has by far the most complex structure, in which case a significant reduction in iteration numbers may be expected without a large additional computational workload per iteration.

\section{Numerical Comparison with Block Diagonal Preconditioner}\label{sec:3}

\begin{table} 
\centering
\caption{Known theoretical results for minimum and maximum (negative and positive) eigenvalues of $\mathcal{P}_{D,k}^{-1}\mathcal{A}_k$, for different values of $k$.}\label{TableBounds}
\begin{tabular}{|c||c|c|c|c|c}
\cline{1-5}
$k$ & Min. --ve & Max. --ve & Min. +ve & Max. +ve & Eigenvalue range (to 3 d.p.) \\ \hline
1 & $-2\,\text{cos}(\frac{\pi}{3})$ & $2\,\text{cos}(\frac{3\pi}{5})$ & $2\,\text{cos}(\frac{\pi}{3})$ & $2\,\text{cos}(\frac{\pi}{5})$ & $[-1,-0.618]\cup[1,1.618]$ \\ \cline{1-5}
2 & $-2\,\text{cos}(\frac{\pi}{5})$ & $2\,\text{cos}(\frac{3\pi}{5})$ & $2\,\text{cos}(\frac{3\pi}{7})$ & $2\,\text{cos}(\frac{\pi}{7})$ & $[-1.618,-0.618]\cup[0.445,1.802]$ \\ \cline{1-5}
3 & $-2\,\text{cos}(\frac{\pi}{7})$ & $2\,\text{cos}(\frac{5\pi}{9})$ & $2\,\text{cos}(\frac{3\pi}{7})$ & $2\,\text{cos}(\frac{\pi}{9})$ & $[-1.802,-0.347]\cup[0.445,1.879]$ \\ \cline{1-5}
\end{tabular}
\end{table}

\begin{figure}
\centering
\begin{subfigure}[b]{0.47\textwidth}
\centering
\includegraphics[width=\textwidth]{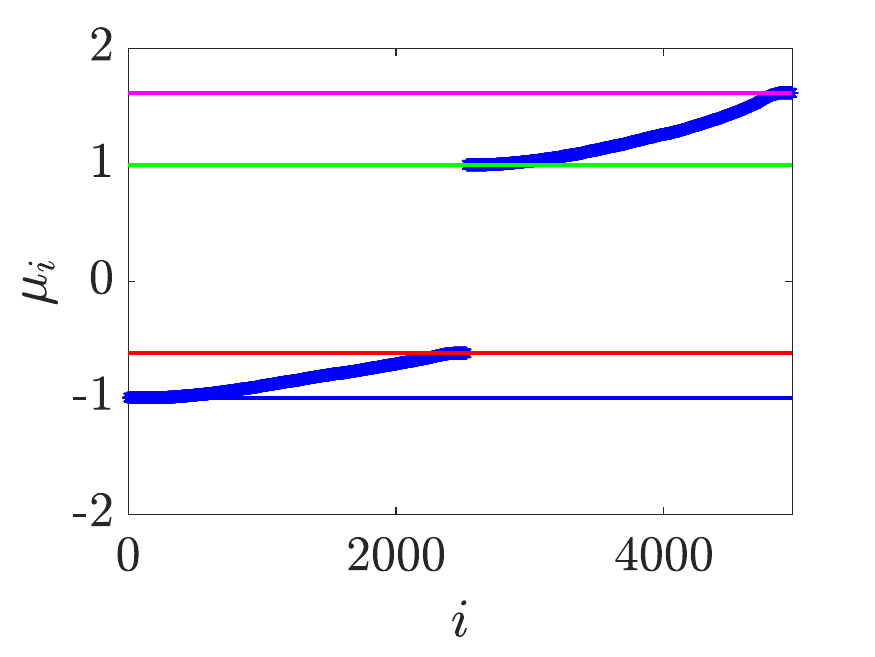}
\end{subfigure}
\begin{subfigure}[b]{0.47\textwidth}
\centering
\includegraphics[width=\textwidth]{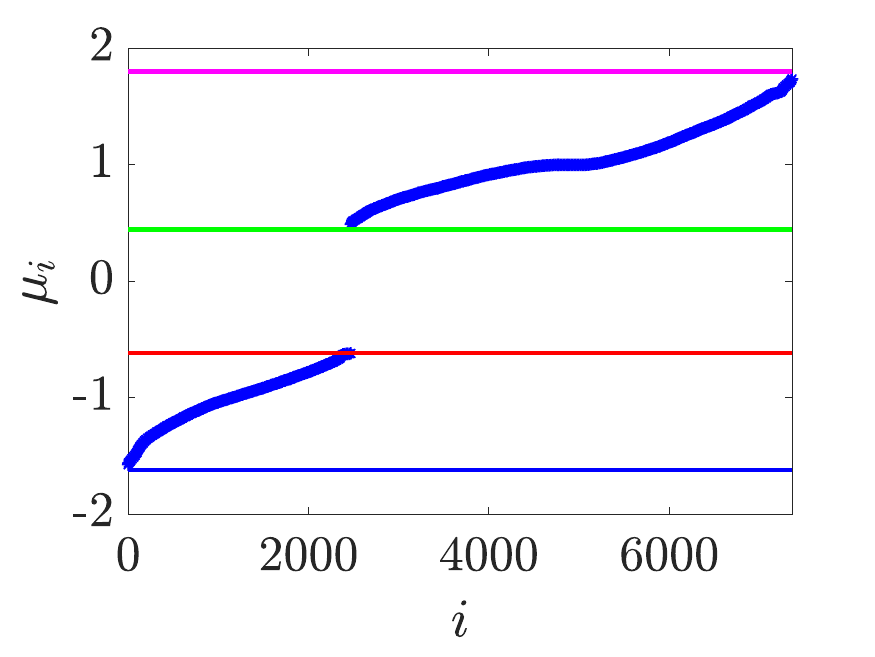}
\end{subfigure}
\begin{subfigure}[b]{0.47\textwidth}
\centering
\includegraphics[width=\textwidth]{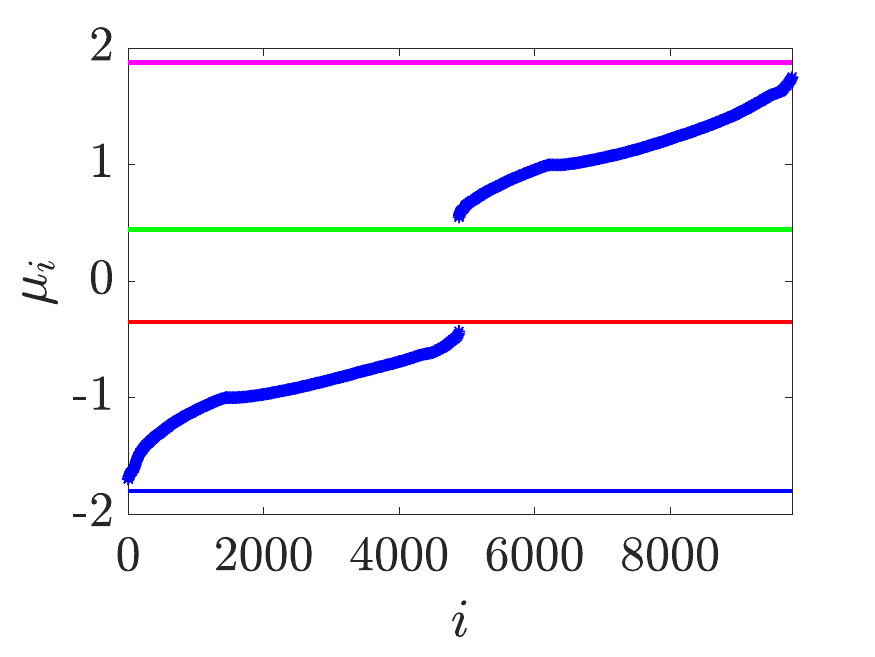}
\end{subfigure}
\caption{Eigenvalues of $\mathcal{P}_{D,k}^{-1}\mathcal{A}_k$ for randomly-generated examples, for $k=1$ \emph{(top-left)}, $k=2$ \emph{(top-right)}, and $k=3$ \emph{(bottom)}. Horizontal lines correspond to theoretical bounds.}
\label{Fig1}
\end{figure}

In order to offer a comparison between the preconditioners $\mathcal{P}_{D,k}$ and $\mathcal{P}_k$, we first recapitulate some known results for $\mathcal{P}_{D,k}$. For $k=1$, that is classical (generalized) saddle-point systems, under the Main Assumption of Section \ref{sec:1}, the eigenvalues of $\mathcal{P}_{D,k}^{-1}\mathcal{A}_k$ are contained within $[-1,\frac{1}{2}(1-\sqrt{5})]\cup[1,\frac{1}{2}(1+\sqrt{5})]$ (see \cite[Corollary 1]{AxNe06}, \cite[Theorem 4]{PearsonDPhil}, and \cite[Lemma 2.2]{SiWa94}). Further, in \cite[Theorem 4]{BrGr21} and \cite[Theorem 5.3]{PearsonPotschka} it is shown that for $k=2$, the bounds are $[-\frac{1}{2}(1+\sqrt{5}),\frac{1}{2}(1-\sqrt{5})]\cup[2\,\text{cos}(\frac{3\pi}{7}),2\,\text{cos}(\frac{\pi}{7})]$. Below, we extend this analysis to the case $k=3$, with a proof given in Appendix \ref{app:proof}:
\begin{Thm}\label{thm:k=3}
Let $k=3$, and assume that the Main Assumption of Section \ref{sec:1} holds. Then all eigenvalues $\mu$ of $\mathcal{P}_{D,k}^{-1}\mathcal{A}_k$ satisfy:
\begin{equation*}
\ \mu\left(\mathcal{P}_{D,k}^{-1}\mathcal{A}_k\right)\in\left[-2\,\emph{cos}\left(\frac{\pi}{7}\right),2\,\emph{cos}\left(\frac{5\pi}{9}\right)\right]\cup\left[2\,\emph{cos}\left(\frac{3\pi}{7}\right),2\,\emph{cos}\left(\frac{\pi}{9}\right)\right],
\end{equation*}
which to $3$ decimal places are $[-1.801,-0.347]\cup[0.445,1.879]$.
\end{Thm}

\begin{figure}
\centering
\begin{subfigure}[b]{0.47\textwidth}
\centering
\includegraphics[width=\textwidth]{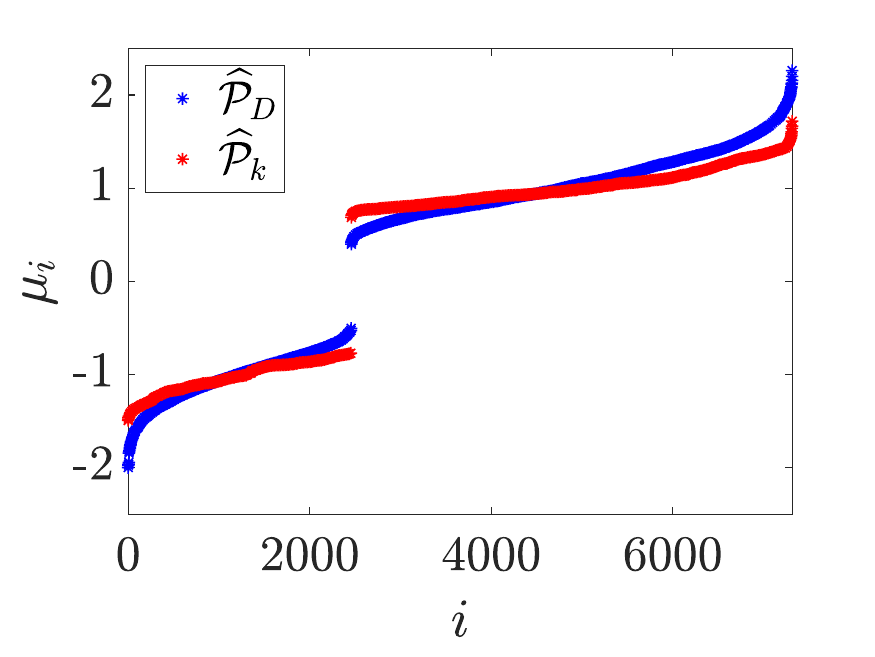}
\end{subfigure}
\begin{subfigure}[b]{0.47\textwidth}
\centering
\includegraphics[width=\textwidth]{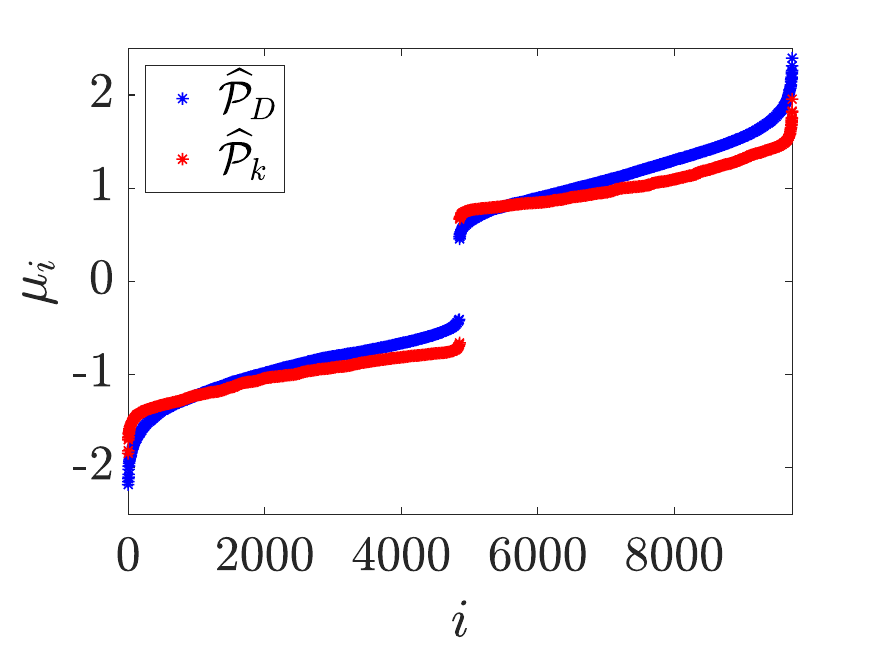}
\end{subfigure}
\begin{subfigure}[b]{0.47\textwidth}
\centering
\includegraphics[width=\textwidth]{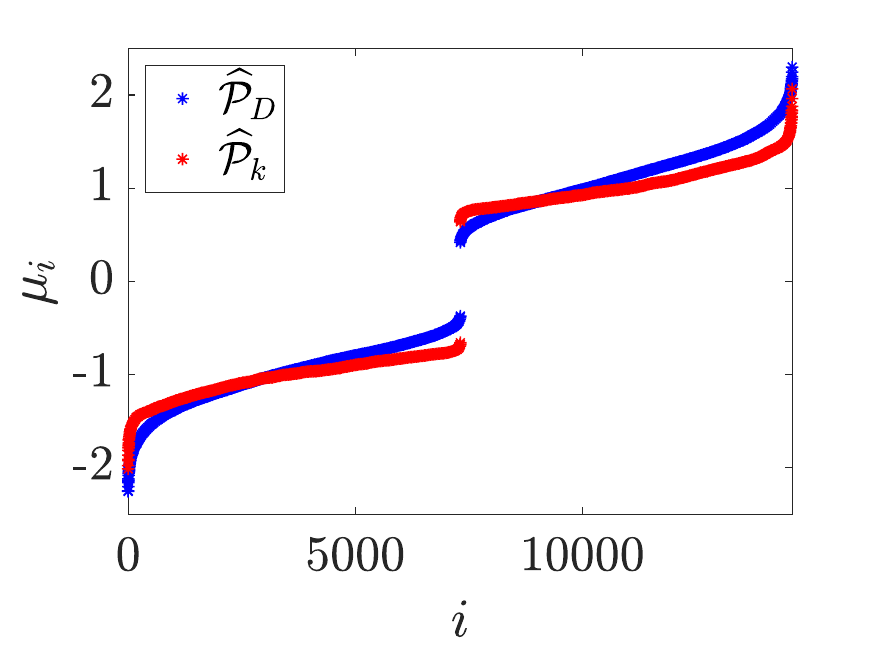}
\end{subfigure}
\begin{subfigure}[b]{0.47\textwidth}
\centering
\includegraphics[width=\textwidth]{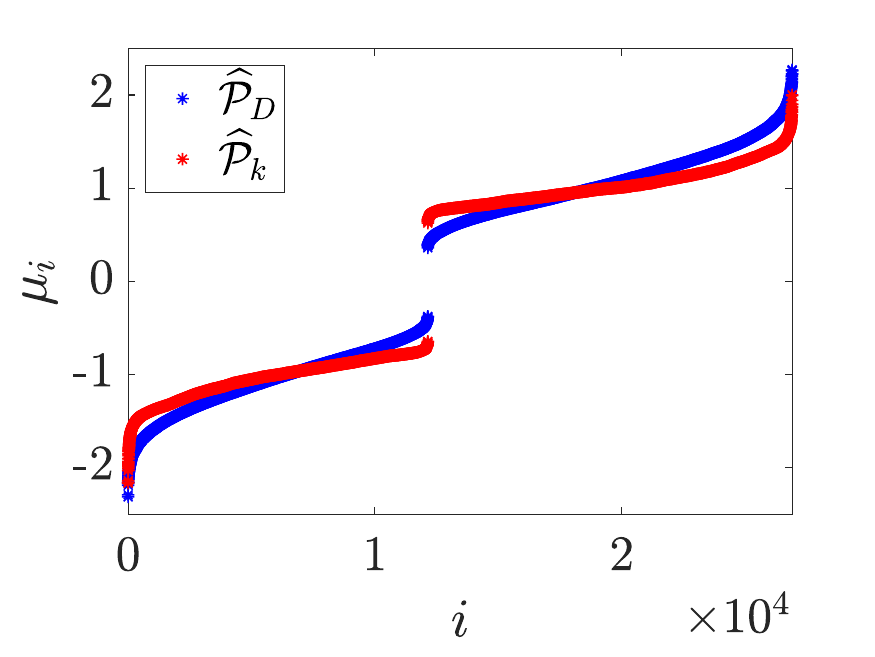}
\end{subfigure}
\caption{Eigenvalues of $\widehat{\mathcal{P}}_{D,k}^{-1}\mathcal{A}_k$ \emph{(blue)} and $\widehat{\mathcal{P}}_k^{-1}\mathcal{A}_k$ \emph{(red)} for randomly-generated examples, for $k=2$ \emph{(top-left)}, $k=3$ \emph{(top-right)}, $k=5$ \emph{(bottom-left)}, and $k=10$ \emph{(bottom-right)}.}
\label{Fig2}
\end{figure}

A summary of the above results are contained in Table \ref{TableBounds}, noting that $1=2\,\text{cos}(\frac{\pi}{3})$, $\frac{1}{2}(1+\sqrt{5})=2\,\text{cos}(\frac{\pi}{5})$, $\frac{1}{2}(1-\sqrt{5})=2\,\text{cos}(\frac{3\pi}{5})$. Some noticeable patterns emerge in the known bounds for varying $k$, and it is likely these could be generalized to higher values of $k$ with even more technical analysis. An important observation is that the eigenvalue bounds appear to spread out, in particular with the eigenvalues of smallest magnitude moving closer to the origin, as $k$ increases. A result of the form of Theorem \ref{thm:Multiple} is therefore highly desirable when one considers the utility of the new preconditioner $\mathcal{P}_k$ for a range of $k$. It is worth pointing out that the majority of the applications considered to date are of the form $k=2$ (see \cite{AliBeikBenzi,BSZ,BrGr21,CJL22,MNN,PearsonPotschka}), with interesting examples with $k=3$ and $k=4$ (in the latter case not block tridiagonal) recently considered in \cite{BSZ}. As this is a rapidly-developing field, we believe the theoretical and practical consideration of multiple saddle-point systems with higher numbers of blocks, and in particular the development of robust preconditioners, will enable and drive forward the consideration of more complex such problems.

\begin{table} 
\centering
\caption{Average {\scshape Minres} iterations required for convergence, with randomly-generated examples, for a range of $k$.}\label{Table0}
\begin{tabular}{|c||c|c|c|c|c|c|c|c|}
\hline
$k$ & 1 & 2 & 3 & 4 & 5 & 10 & 15 & 20 \\ \hline
Average DoF & 493 & 754 & 1,007 & 1,247 & 1,498 & 2,738 & 3,977 & 5,240 \\ \hline \hline
$\widehat{\mathcal{P}}_{D,k}$ iterations & 33.1 & 59.9 & 65.6 & 74.1 & 74.1 & 80.4 & 80.0 & 80.8 \\ \hline
$\widehat{\mathcal{P}}_k$ iterations & 30.4 & 34.0 & 35.0 & 34.6 & 34.8 & 34.3 & 33.6 & 33.6 \\ \hline
\end{tabular}
\end{table}

In Figure \ref{Fig1}, we show how the eigenvalues of $\mathcal{P}_{D,k}^{-1}\mathcal{A}_k$ behave in practice for a sequence of randomly-constructed examples, for $k=1,2,3$. In more detail, for each value of $k$ we compute 100 random matrices in {\scshape Matlab} R2018a, and plot all eigenvalues obtained. We select random dimensions $n_j$ using \texttt{fix(20+10*rand)}. To construct matrices $A_j$, we take the symmetric part of random matrices (formed using \texttt{randn}), then add the absolute value of the smallest negative eigenvalue multiplied by an identity matrix, except for $A_0$ in which case we add 1.01 times this matrix to ensure positive definiteness. We also construct $B_j$ using the \texttt{randn} function. The horizontal lines in Figure \ref{Fig1} correspond to the theoretical bounds in Table \ref{TableBounds}.

We also examine the behaviour of the eigenvalues of the preconditioned systems when $A_0$ and the Schur complements are applied inexactly within preconditioners $\widehat{\mathcal{P}}_{D,k}$ and $\widehat{\mathcal{P}}_k$. To present Figure \ref{Fig2}, for a number of values of $k$ we construct 100 random test problems as above. Letting $\mu_{\max}$ and $\mu_{\min}$ denote the maximum and minimum eigenvalues of $A_0$, we set $\widehat{S}_0:=\widehat{A}_0=\frac{1}{\mu_{\max}-\mu_{\min}}[(\frac{2}{3}\mu_{\max}-2\mu_{\min})A_0+\frac{4}{3}\mu_{\max}\mu_{\min}I_{n_0}]$, guaranteeing that the eigenvalues of $\widehat{A}_0^{-1}A_0$ are contained in $[\frac{1}{2},\frac{3}{2}]$; we then set $\widehat{S}_j=A_j+B_j \widehat{S}_{j-1}^{-1} B_j^T$ for $j=1,...,k$. This tests the effect of inexactness in the approximation of $A_0$ `propagating' through the preconditioners $\mathcal{P}_{D,k}$ and $\mathcal{P}_k$. Figure \ref{Fig2} shows that one may frequently expect more of the eigenvalues of $\widehat{\mathcal{P}}_k^{-1}\mathcal{A}_k$ to be clustered around $1$ and $-1$ than those of $\widehat{\mathcal{P}}_{D,k}^{-1}\mathcal{A}_k$, in particular with more small-magnitude eigenvalues further from the origin. Table \ref{Table0} presents the average number of {\scshape Minres} iterations, using the routine at \cite{minres_code}, required to solve 100 (different) randomly-generated test problems with $\widehat{\mathcal{P}}_{D,k}$ and $\widehat{\mathcal{P}}_k$, for a range of $k$ and with $n_j$ now defined using \texttt{fix(200+100*rand)}. We solve to a tolerance of $10^{-10}$, and refer to \cite{FoSa12} for a justification of the stopping condition used -- this criterion is motivated using backward error analysis, and measures whether the 2-norm of the current residual is as small as the specified tolerance, multiplied by a factor involving an estimation of $\lVert\mathcal{A}\rVert_2 \lVert\widehat{x}\rVert_2$, where $\widehat{x}$ denotes the approximate solution at the current iterate (see also \cite{minres_code}). The same stopping criterion is used for both preconditioners to ensure a fair comparison. We also report the average number of degrees of freedom (DoF's) for each $k$. Although this is a synthetic test case, as one would not exactly factor and then apply Schur complement approximations in practice, the significantly smaller numbers of iterations routinely required by $\widehat{\mathcal{P}}_k$ (typically less than one half of those required by $\widehat{\mathcal{P}}_{D,k}$) further indicate the potential utility of the new preconditioner.

\section{Applications and Numerical Results}\label{sec:4}

As further proofs of concept of our new preconditioning strategy, we now present two applications of multiple saddle-point systems from PDE-constrained optimization, and solve them numerically using {\scshape Minres} with both preconditioners $\mathcal{P}_{D,k}$ and $\mathcal{P}_k$, using the same approximations for $A_0$ and the Schur complements. All matrices are generated in Python using the package \cite{skfem2020}, whereupon we solve the resulting systems in {\scshape Matlab} R2018a using the {\scshape Minres} routine at \cite{minres_code} to tolerance $10^{-10}$. All tests are carried out on an Intel(R) Core(TM) i7-6700T CPU @ 2.80GHz quad-core processor, and we report the DoF's of the systems solved.

\subsection{Double saddle-point system from PDE-constrained optimization problem}\label{sec:4-1}

To illustrate the potential improvements of the preconditioner $\mathcal{P}_k$ proposed in this paper as opposed to the block diagonal preconditioner $\mathcal{P}_{D,k}$, even in the case $k=2$, we consider the discretization of the following PDE-constrained optimization problem with boundary observations:
\begin{align*}
\min_{u,f}~~&\frac{1}{2}\left\|u-\widehat{u}\right\|_{L^2(\partial\Omega)}^2+\frac{\alpha}{2}\left\|f\right\|_{L^2(\Omega)}^2 \\
\text{s.t.}~~~&\left\{\begin{array}{rl}
-\Delta u+u+f=0 & \text{in }\Omega, \\
\frac{\partial u}{\partial n}=0 & \text{on }\partial\Omega, \\
\end{array}\right.
\end{align*}
posed on a domain $\Omega$ with boundary $\partial\Omega$, with $u$ and $f$ the \emph{state} and \emph{control variables}, $\widehat{u}$ a given \emph{desired state}, and $\alpha>0$ a given regularization parameter. This problem was examined by the authors of \cite{MNN}, and problems of similar structure were also considered in \cite{SognZulehner}. As in \cite{MNN}, we discretize this problem using $H^1$-conforming piecewise linear Lagrange finite elements $\phi_i$ on a triangularized mesh, leading to a linear system of the form
\begin{equation}\label{PDECO_system}
\left(\begin{array}{ccc}
\alpha M_h & M_h & O \\ M_h & O & L_h \\ O & L_h & Q_h \\
\end{array}\right)\left(\begin{array}{c}
f_h \\ p_h \\ u_h \\
\end{array}\right)=\left(\begin{array}{c}
0 \\ 0 \\ \widehat{u}_h \\
\end{array}\right),
\end{equation}
with $M_h \, (\, =M_h^T \,)$ a finite element \emph{mass matrix}, $L_h \,  (\, =L_h^T \,)$ the sum of a \emph{stiffness matrix} and a mass matrix, and $Q_h$ a boundary mass matrix. Here, $u_h$, $f_h$ are the discretized versions of $u$, $f$, with $p_h$ the corresponding discretized \emph{adjoint variable}, and $\widehat{u}_h$ contains terms of the form $\int_{\partial\Omega}\widehat{u}\phi_i$. For our numerical tests, we consider the above problem on the domain $\Omega=(0,1)^2$, with $\widehat{u}$ generated by solving the forward PDE with `true' control $4x(1-x)+y$, $x$ and $y$ denoting the spatial variables, for a range of mesh parameters $h$ and values of $\alpha$.

\begin{table} 
\centering
\caption{Double saddle-point problem: {\scshape Minres} iterations (top of each cell) and CPU time in seconds (bottom) required for convergence with preconditioner $\mathcal{P}_{D,k}$.}\label{Table1}
\begin{tabular}{|c|c||c|c|c|c|c|}
\hline
DoF & $h\backslash\alpha$ & $1$ & $10^{-1}$ & $10^{-2}$ & $10^{-3}$ & $10^{-4}$ \\ \hline \hline
\multirow{2}{*}{867} & \multirow{2}{*}{$2^{-4}$} & 17 & 21 & 24 & 27 & 20 \\
 & & 0.0109 & 0.0137 & 0.0165 & 0.0228 & 0.0156 \\ \hline
\multirow{2}{*}{3,267} & \multirow{2}{*}{$2^{-5}$} & 17 & 21 & 22 & 26 & 18 \\
 & & 0.0397 & 0.0437 & 0.0446 & 0.0515 & 0.0361 \\ \hline
\multirow{2}{*}{12,675} & \multirow{2}{*}{$2^{-6}$} & 14 & 19 & 22 & 25 & 15 \\
 & & 0.118 & 0.153 & 0.172 & 0.193 & 0.118 \\ \hline
\multirow{2}{*}{49,923} & \multirow{2}{*}{$2^{-7}$} & 14 & 19 & 21 & 20 & 14 \\
 & & 0.463 & 0.618 & 0.676 & 0.619 & 0.448 \\ \hline
\multirow{2}{*}{198,147} & \multirow{2}{*}{$2^{-8}$} & 16 & 18 & 21 & 17 & 12 \\
 & & 2.04 & 2.29 & 2.62 & 2.16 & 1.54 \\ \hline
 \multirow{2}{*}{789,507} & \multirow{2}{*}{$2^{-9}$} & 14 & 18 & 19 & 14 & 12 \\
 & & 7.43 & 9.26 & 9.74 & 7.38 & 6.44 \\ \hline
 \multirow{2}{*}{3,151,875} & \multirow{2}{*}{$2^{-10}$} & 14 & 17 & 18 & 15 & 11 \\
 & & 31.6 & 37.7 & 39.7 & 33.6 & 25.8 \\ \hline
\end{tabular}
\end{table}

\begin{table} 
\centering
\caption{Double saddle-point problem: {\scshape Minres} iterations (top of each cell) and CPU time in seconds (bottom) required for convergence with preconditioner $\mathcal{P}_k$.}\label{Table2}
\begin{tabular}{|c|c||c|c|c|c|c|}
\hline
DoF & $h\backslash\alpha$ & $1$ & $10^{-1}$ & $10^{-2}$ & $10^{-3}$ & $10^{-4}$ \\ \hline \hline
\multirow{2}{*}{867} & \multirow{2}{*}{$2^{-4}$} & 8 & 9 & 11 & 12 & 12 \\
 & & 0.0106 & 0.00900 & 0.00889 & 0.0124 & 0.0106 \\ \hline
\multirow{2}{*}{3,267} & \multirow{2}{*}{$2^{-5}$} & 8 & 9 & 9 & 12 & 9 \\
 & & 0.0208 & 0.0216 & 0.0216 & 0.0285 & 0.0210 \\ \hline
\multirow{2}{*}{12,675} & \multirow{2}{*}{$2^{-6}$} & 7 & 9 & 9 & 12 & 8 \\
 & & 0.0658 & 0.0836 & 0.0830 & 0.105 & 0.0728 \\ \hline
\multirow{2}{*}{49,923} & \multirow{2}{*}{$2^{-7}$} & 7 & 9 & 9 & 10 & 7 \\
 & & 0.290 & 0.357 & 0.352 & 0.363 & 0.267 \\ \hline
\multirow{2}{*}{198,147} & \multirow{2}{*}{$2^{-8}$} & 7 & 7 & 9 & 10 & 7 \\
 & & 1.11 & 1.10 & 1.36 & 1.46 & 1.08 \\ \hline
 \multirow{2}{*}{789,507} & \multirow{2}{*}{$2^{-9}$} & 7 & 7 & 9 & 9 & 7 \\
 & & 4.50 & 4.50 & 5.55 & 5.55 & 4.50 \\ \hline
 \multirow{2}{*}{3,151,875} & \multirow{2}{*}{$2^{-10}$} & 7 & 7 & 8 & 7 & 7 \\
 & & 19.8 & 19.7 & 22.1 & 19.7 & 19.7 \\ \hline
\end{tabular}
\end{table}

\begin{table} 
\centering
\caption{Double saddle-point problem: {\scshape Minres} iterations required for convergence with preconditioners $\mathcal{P}_{D,k}$, $\mathcal{P}_k$, for varying numbers of Chebyshev semi-iterations.}\label{Table3}
\begin{tabular}{|c||c|c|c|c|c|c|c|c||c|c|c|c|c|c|c|c|}
\hline
 & \multicolumn{8}{c||}{$\mathcal{P}_{D,k}$, \# Chebyshev iterations =} & \multicolumn{8}{c|}{$\mathcal{P}_k$, \# Chebyshev iterations =} \\ \cline{2-17}
$h$ & 1 & 2 & 3 & 4 & 5 & 7 & 10 & 20 & 1 & 2 & 3 & 4 & 5 & 7 & 10 & 20 \\ \hline \hline
$2^{-4}$ & 42 & 37 & 25 & 23 & 24 & 23 & 23 & 23 & 40 & 19 & 12 & 11 & 11 & 9 & 9 & 8 \\ \hline
$2^{-5}$ & 37 & 34 & 23 & 23 & 22 & 22 & 22 & 22 & 36 & 17 & 13 & 11 & 9 & 9 & 9 & 8 \\ \hline
$2^{-6}$ & 34 & 29 & 22 & 22 & 22 & 21 & 21 & 21 & 31 & 16 & 11 & 10 & 9 & 9 & 9 & 8 \\ \hline
$2^{-7}$ & 30 & 28 & 21 & 21 & 21 & 21 & 21 & 21 & 27 & 15 & 11 & 10 & 9 & 9 & 9 & 8 \\ \hline
$2^{-8}$ & 27 & 26 & 21 & 21 & 21 & 21 & 21 & 21 & 24 & 15 & 10 & 9 & 9 & 9 & 9 & 8 \\ \hline
$2^{-9}$ & 25 & 23 & 19 & 19 & 19 & 19 & 19 & 19 & 21 & 13 & 9 & 9 & 9 & 9 & 7 & 7 \\ \hline
$2^{-10}$ & 21 & 24 & 18 & 18 & 18 & 18 & 18 & 18 & 18 & 13 & 8 & 8 & 8 & 8 & 7 & 7 \\ \hline
\end{tabular}
\end{table}

In the notation of multiple saddle-point systems outlined in this paper, the matrix in \eqref{PDECO_system} leads to
\begin{align*}
&A_0=\alpha M_h, \quad A_1=O, \quad A_2=Q_h, \quad B_1=B_1^T=M_h, \quad B_2=B_2^T=L_h, \\
&S_1=\frac{1}{\alpha}M_h, \quad S_2=Q_h+\alpha L_h M_h^{-1} L_h.
\end{align*}

For our first tests, within $\widehat{A}_0$ and $\widehat{S}_1$ we apply 5 Chebyshev semi-iterations (see \cite{GVI,GVII,WathenRees}) with Jacobi splitting to approximate the action of $M_h^{-1}$ on a vector. To approximate $S_2$ we take $\widehat{S}_2=\alpha L_h M_h^{-1} L_h$, and use 2 V-cycles of the \texttt{HSL\_MI20} algebraic multigrid solver \cite{HSL_MI20,HSL_MI20_code} (with 2 symmetric Gauss--Seidel iterations as pre-/post-smoother) to approximate the action of $L_h^{-1}$ on each occasion. We select these numbers of Chebyshev semi-iterations and multigrid V-cycles as we find these compute inner approximations to satisfactory accuracy, while increasing these numbers only mildly affects the outer {\scshape Minres} iterations (see Table \ref{Table3} for evidence of this for Chebyshev semi-iterations, which we note are cheap to compute compared to multigrid V-cycles). In Table \ref{Table1}, we report the numbers of preconditioned {\scshape Minres} iterations required for convergence with the (approximately applied) block diagonal preconditioner $\mathcal{P}_{D,k}$, as well as CPU times (reported in seconds), for a range of $h$ and $\alpha$.\footnote{We note that the 2-norm of the residual, scaled by the square root of the dimension of the problem, was smaller than $10^{-4}$ upon convergence for all choices of $h$ and $\alpha$, with either choice of preconditioner; this measure decreased to $\mathcal{O}(10^{-5})$ for some problems with the smallest $\alpha=10^{-4}$, down to $\mathcal{O}(10^{-8})$ with the largest $\alpha=1$. For 25 out of the 35 choices of $h$ and $\alpha$, a smaller residual norm was obtained at the point of convergence when using $\mathcal{P}_k$ than with $\mathcal{P}_{D,k}$.} In Table \ref{Table2}, we report the same quantities with the (approximately applied) preconditioner $\mathcal{P}_k$. We observe that both the iteration numbers and CPU times are lower in each case with preconditioner $\mathcal{P}_k$, often by a factor of roughly two. This demonstrates the efficacy of the preconditioning approach outlined in this paper. Of course, there are many choices of approximations for the individual blocks, and the relative merits of $\mathcal{P}_{D,k}$ and $\mathcal{P}_k$ will of course depend on these precise choices.


As emphasized above, we anticipate the benefit of $\mathcal{P}_k$ to be more obvious when $A_0$, $S_1$ may be well approximated in a computationally cheap way, in comparison to $S_2$, as any inexactness in these approximations is approximated twice within $\mathcal{P}_k$ but only once within $\mathcal{P}_{D,k}$. Indeed we observe this feature for problem \eqref{PDECO_system}: in Table \ref{Table3} we present iteration numbers using both preconditioners for a fixed $\alpha=10^{-2}$, while varying the number of Chebyshev semi-iterations used to approximate $M_h$, to test the effect of deliberately taking a worse approximation of $A_0$, $S_1$. We find that $\mathcal{P}_{D,k}$ converges in roughly the same number of iterations as $\mathcal{P}_k$ when only one iteration of Chebyshev semi-iteration is applied to $M_h$; this is equivalent to a scaled diagonal approximation of $M_h$, which is mesh-robust \cite{WathenEigs} but relatively inaccurate. When two or more Chebyshev semi-iterations are applied to $M_h$, {\scshape Minres} converges in both fewer iterations and lower CPU time when accelerated by $\mathcal{P}_k$ than by $\mathcal{P}_{D,k}$. We highlight that the case of taking a large number of Chebyshev semi-iterations resembles the setting of Theorem \ref{thm:ApproxSk}, in the sense that the blocks $A_0$, $S_1$ may be approximated to arbitrary accuracy, but this cannot be done for the final Schur complement $S_2$ in a computationally cheap way.

\subsection{Quadruple saddle-point system from PDE-constrained optimization with state constraints}

The presence of state constraints in PDE-constrained optimization poses severe difficulties on the analytical as well as on the numerical level (see, e.g., \cite{ItKu2003}), because in general the Lagrange multipliers are measures and thus exhibit only low regularity. We consider the problem
\begin{align*}
  \min_{u \in H^1(\Omega), f \in L^2(\Omega)}~~&\frac{1}{2}\left\|u-\widehat{u}\right\|_{L^2(\Omega)}^2+\frac{\alpha}{2}\left\|f\right\|_{L^2(\Omega)}^2 \\
  \text{s.t.}~~~&\left\{\begin{array}{rl}
  -\Delta u = f & \text{in }\Omega, \\
  \frac{\partial u}{\partial n}=0 & \text{on }\partial\Omega, \\
  u \le \bar{u} & \text{a.e. in } \Omega, \\
  \end{array}\right.
\end{align*}
with a continuous upper state bound $\bar{u} \in C(\overline{\Omega})$. Instead of using a regularization approach as in \cite{ItKu2003}, we introduce to the original problem an additional ``lifting'' variable $\widetilde{u} \in L^2(\Omega)$ according to
\begin{align*}
  \min_{u \in H^1(\Omega), \widetilde{u} \in L^2(\Omega), f \in L^2(\Omega)}~~&\frac{1}{2}\left\|u-\widehat{u}\right\|_{L^2(\Omega)}^2+\frac{\alpha}{2}\left\|f\right\|_{L^2(\Omega)}^2 \\
  \text{s.t.}~~~&\left\{\begin{array}{rll}
  -\Delta u = f & \text{in }\Omega, & | \cdot p\\
  \frac{\partial u}{\partial n}=0 & \text{on }\partial\Omega, \\
  u = \widetilde{u} & \text{in }\Omega, & | \cdot \widetilde{p} \\
  \widetilde{u} \le \bar{u} & \text{a.e. in } \Omega, \\
  \end{array}\right.
\end{align*}
where the last two constraints are supposed to hold in the larger space $L^2(\Omega)$, so that a pointwise interpretation of necessary optimality conditions becomes possible. The corresponding adjoint variables to the equality constraints are $p \in H^1(\Omega)$ and $\widetilde{p} \in L^2(\Omega)$. We can then apply the sequential homotopy method of \cite{PoBo21} with the limitation that due to a norm gap for the pointwise defined projection operator (see, e.g., \cite[Section 3.3]{Ulbrich2011}), semismoothness of the resulting non-smooth systems can only be expected for each discretized problem but not in function space. Hence, the overall nonlinear, nonsmooth outer iterations are mesh-dependent (see Table \ref{TableOuterIterations}). However, the required number of outer homotopy iterations on the finest considered mesh are moderate, hence fast methods for the arising linear subsystems are of interest.

\begin{table}
  \centering
  \caption{Quadruple saddle-point problem: The outer sequential homotopy iterations are mesh-dependent. Computations were performed for a termination tolerance of $10^{-10}$ with $\lambda = 10^{-10}$ and direct solution of the linear subproblems.}
  \begin{tabular}{|c||c|c|c|c|c|}
    \hline
    $h$ & $2^{-4}$ & $2^{-5}$ & $2^{-6}$ & $2^{-7}$ & $2^{-8}$ \\
    \hline
    Seq.~Hom.~Iter. & 4 & 4 & 7 & 13 & 26 \\
    \hline
  \end{tabular}
  \label{TableOuterIterations}
\end{table}

\begin{table} 
\centering
\caption{Quadruple saddle-point problem: {\scshape Minres} iterations (top of each cell) and CPU time in seconds (bottom) required for convergence with preconditioner $\mathcal{P}_{D,k}$.}\label{Table4}
\begin{tabular}{|c|c||c|c|c||c|c|c|}
\hline
 &  & \multicolumn{3}{c||}{$\lambda=10^{-7}$, $\alpha=$} & \multicolumn{3}{c|}{$\lambda=5\times10^{-9}$, $\alpha=$} \\ \cline{3-8}
DoF & $h$ & $10^{-6}$ & $10^{-8}$ & $10^{-10}$ & $10^{-6}$ & $10^{-8}$ & $10^{-10}$ \\ \hline \hline
\multirow{2}{*}{9,852} & \multirow{2}{*}{$2^{-5}$} & 70 & 60 & 60 & 70 & 77$^\dag$ & 48$^\dag$ \\
 & & 0.553 & 0.506 & 0.456 & 0.601 & 0.150 & 0.0951 \\ \hline
\multirow{2}{*}{38,832} & \multirow{2}{*}{$2^{-6}$} & 72 & 60 & 60 & 72 & 52 & 46 \\
 & & 2.18 & 1.82 & 1.85 & 2.13 & 1.65 & 1.33 \\ \hline
\multirow{2}{*}{154,224} & \multirow{2}{*}{$2^{-7}$} & 85 & 62 & 62 & 85 & 51 & 46 \\
 & & 10.1 & 7.10 & 7.22 & 9.63 & 6.15 & 5.79 \\ \hline
\multirow{2}{*}{614,664} & \multirow{2}{*}{$2^{-8}$} & 110 & 79 & 77 & 109 & 57 & 48 \\
 & & 49.9 & 36.6 & 35.3 & 50.6 & 26.4 & 23.1 \\ \hline
\end{tabular}
\end{table}

With the variable ordering $(f, p, u, \widetilde{p}, \widetilde{u})$, the discretized linear subproblems (after elimination of active variables and nonlinear transformation to avoid the augmented Lagrangian matrix terms) have the form
\begin{equation}\label{PDECO_quadruple}
\left(\begin{array}{ccccc}
(\alpha+\lambda)M_h & -M_h & O & O & O \\ -M_h & -\frac{\lambda}{1+\rho\lambda}L_h & K_h & O & O \\ O & K_h & M_h+\lambda L_h & -M_h & O \\ O & O & -M_h & -\frac{\lambda}{1+\rho\lambda}M_h & (M_h^{(i,:)})^T \\ O & O & O & M_h^{(i,:)} & \lambda M_h^{(i,i)} \\
\end{array}\right)
\begin{pmatrix}
  f_h\\ p_h\\ u_h\\ \widetilde{p}_h\\ \widetilde{u}^{(i)}_h
\end{pmatrix}
=b_h,
\end{equation}
which is a multiple saddle-point system with $k=4$. Here, $K_h$ is a stiffness matrix, $M_h^{(i,i)}$ is a mass matrix defined over the inactive nodes, and $M_h^{(i,:)}$ is the analogous matrix with rows corresponding to inactive nodes. We also denote $K_h^{(i,i)}$ as a stiffness matrix defined over the inactive nodes. 
The right-hand side is computed by finite element assembly and is nonzero except for the third block component, in which $\widehat{u}$ enters. 
The parameter $\lambda>0$ arises within the sequential homotopy method as the reciprocal $\lambda = \frac{1}{\Delta t}$ of a projected backward Euler timestep of size $\Delta t > 0$ on a projected gradient/antigradient flow on the augmented Lagrangian, and $\rho>0$ is the standard augmented Lagrangian parameter. It is set here to a relatively low value of $\rho=10^{-5}$ for these experiments because this ensured rapid convergence of the outer iteration. Other problems might require higher values of $\rho$, and both preconditioners also perform well for these values.

\begin{table} 
\centering
\caption{Quadruple saddle-point problem: {\scshape Minres} iterations (top of each cell) and CPU time in seconds (bottom) required for convergence with preconditioner $\mathcal{P}_k$.}\label{Table5}
\begin{tabular}{|c|c||c|c|c||c|c|c|}
\hline
 &  & \multicolumn{3}{c||}{$\lambda=10^{-7}$, $\alpha=$} & \multicolumn{3}{c|}{$\lambda=5\times10^{-9}$, $\alpha=$} \\ \cline{3-8}
DoF & $h$ & $10^{-6}$ & $10^{-8}$ & $10^{-10}$ & $10^{-6}$ & $10^{-8}$ & $10^{-10}$ \\ \hline \hline
\multirow{2}{*}{9,852} & \multirow{2}{*}{$2^{-5}$} & 21 & 19 & 19 & 21 & 96$^\dag$ & 40$^\dag$ \\
 & & 0.299 & 0.258 & 0.252 & 0.306 & 0.275 & 0.116 \\ \hline
\multirow{2}{*}{38,832} & \multirow{2}{*}{$2^{-6}$} & 27 & 20 & 20 & 27 & 17 & 17 \\
 & & 1.36 & 1.04 & 1.06 & 1.33 & 0.921 & 0.867 \\ \hline
\multirow{2}{*}{154,224} & \multirow{2}{*}{$2^{-7}$} & 34 & 25 & 25 & 34 & 18 & 16 \\
 & & 6.56 & 4.93 & 4.96 & 6.45 &  3.74 & 3.48 \\ \hline
\multirow{2}{*}{614,664} & \multirow{2}{*}{$2^{-8}$} & 54 & 31 & 31 & 46 & 24 & 20 \\
 & & 41.2 & 24.2 & 24.6 & 35.4 & 19.0 & 16.5 \\ \hline
\end{tabular}
\end{table}

In the notation of this paper, the matrix \eqref{PDECO_quadruple} gives
\begin{align*}
&A_0=(\alpha+\lambda)M_h, \quad A_1=\frac{\lambda}{1+\rho\lambda}L_h, \quad A_2=M_h+\lambda L_h, \quad A_3=\frac{\lambda}{1+\rho\lambda}M_h, \\
\ &A_4=\lambda M_h^{(i,i)}, \quad B_1=B_1^T=B_3=B_3^T=-M_h, \quad B_2=B_2^T=K_h, \quad B_4=M_h^{(i,:)}, \\
&S_1=\frac{\lambda}{1+\rho\lambda}L_h+\frac{1}{\alpha+\lambda}M_h, \quad S_2=M_h+\lambda L_h+K_h\left[\frac{\lambda}{1+\rho\lambda}L_h+\frac{1}{\alpha+\lambda}M_h\right]^{-1}K_h, \\
&S_3=\frac{\lambda}{1+\rho\lambda}M_h+M_h S_2^{-1} M_h, \quad S_4=\lambda M_h^{(i,i)}+M_h^{(i,:)} S_3^{-1} (M_h^{(i,:)})^T.
\end{align*}
This is a highly complicated system, for which we build preconditioners based on approximations tailored to problems with small values of $\alpha$ and $\lambda$, a regime of significant numerical interest. The key objective here is to compare the preconditioners $\mathcal{P}_{D,k}$ and $\mathcal{P}_k$. Our approximation $\widehat{A}_0$ involves 5 Chebyshev semi-iterations approximating the action of $M_h^{-1}$; $S_1$ is approximated by $\frac{1}{\alpha+\lambda}M_h$, which is then handled in the same way. We approximate $S_2$ by neglecting terms multiplied by $\lambda$ and applying the `matching strategy' of \cite{PW12,PW13} to approximate:
\begin{equation*}
S_2 \approx \underbrace{M_h+(\alpha+\lambda) K_h M_h^{-1} K_h}_{\widetilde{S}_2} \approx \frac{1}{\sqrt{2}}\left(M_h+\sqrt{\alpha+\lambda} \, K_h\right)M_h^{-1}\left(M_h+\sqrt{\alpha+\lambda} \, K_h\right)=:\widehat{S}_2.
\end{equation*}
Using the working of \cite{PW13} we may show that the eigenvalues of $\widehat{S}_2^{-1}\widetilde{S}_2$ are contained within $[\frac{1}{\sqrt{2}},\sqrt{2}]$. The factor of $\frac{1}{\sqrt{2}}$ within $\widehat{S}_2$ is present to give symmetry about $1$ (in a geometric sense) of these eigenvalue bounds. To apply $\widehat{S}_2^{-1}$ we use 2 V-cycles of the \texttt{HSL\_MI20} multigrid solver to approximate $M_h+\sqrt{\alpha+\lambda} \, K_h$ in both instances. To approximate the action of $S_3^{-1}$, we neglect the first term in $S_3$ and write $\widehat{S}_3^{-1}=M_h^{-1} [M_h+(\alpha+\lambda)K_h M_h^{-1} K_h] M_h^{-1}$, where $M_h^{-1}$ is then approximated using 5 Chebyshev semi-iterations on three occasions. Finally, we approximate
\begin{align*}
S_4 \approx{}& M_h^{(i,:)} S_3^{-1} (M_h^{(i,:)})^T \approx M_h^{(i,:)} M_h^{-1} S_2 M_h^{-1} (M_h^{(i,:)})^T \\
\approx{}& \frac{1}{\sqrt{2}}\left(M_h^{(i,i)}+\sqrt{\alpha+\lambda} \, K_h^{(i,i)}\right)\left(M_h^{(i,i)}\right)^{-1}\left(M_h^{(i,i)}+\sqrt{\alpha+\lambda} \, K_h^{(i,i)}\right)=:\widehat{S}_4,
\end{align*}
and use 2 V-cycles of the \texttt{HSL\_MI20} solver to twice approximate $M_h^{(i,i)}+\sqrt{\alpha+\lambda} \, K_h^{(i,i)}$.\footnote{We make an exception to this strategy, indicated by $^\dag$ in Tables \ref{Table4} and \ref{Table5}, when coarsening failed within the \texttt{HSL\_MI20} routine due to the very specific parameter regime, in which case we instead approximate $S_2^{-1}$ and $S_4^{-1}$ with 5 Chebyshev semi-iterations applied to $M_h$ and $M_h^{(i,i)}$, respectively.} We again emphasize that other choices of approximations for the blocks may yield good results.

For the numerical computations, the matrices are obtained using $H^1$-conforming piecewise linear Lagrange finite elements on a triangular mesh approximating the unit disc $\Omega = \{ \mathbf{x} \in \mathbb{R}^2 \mid \lVert \mathbf{x} \rVert_2 < 1 \}$. The problem data are  $\widehat{u}(\mathbf{x}) = 1 - x^2 - y^2$ and $\bar{u}(\mathbf{x}) = \frac{1}{2}$. For easier comparison, we artificially assume for the remaining numerical experiments that the current active set consists of all mesh nodes with distance at most $\frac{1}{2}$ from the origin.

In Tables \ref{Table4} and \ref{Table5} we report the required preconditioned {\scshape Minres} iterations and CPU times with (approximately applied) preconditioners $\mathcal{P}_{D,k}$ and $\mathcal{P}_k$, respectively, as well as DoF's of the problems solved. We note that, once again, our new structure of preconditioner $\mathcal{P}_k$ noticeably outperforms the block diagonal preconditioner $\mathcal{P}_{D,k}$ in terms of iteration numbers. Further, apart from two problem instances with a very coarse discretization, the CPU time is lower when $\mathcal{P}_k$ is used. Indeed for large problems the CPU time is significantly and consistently lower. This further justifies the potential utility of this new preconditioner.




\section{Conclusion}

We proposed a symmetric positive definite preconditioner for block tridiagonal multiple saddle-point systems, such that the preconditioned matrix has only two distinct eigenvalues when the preconditioner is applied exactly, independently of the number of matrix blocks. In the setting where the final Schur complement is the most difficult to cheaply and accurately approximate, which often arises in practice, the approach outlined here can offer a significant advantage over the widely-used block diagonal preconditioner. Valuable future work would include the analysis of the preconditioner when individual blocks are approximated, for different values of $k$, and implementation of the solver for different physical applications.

\section*{Acknowledgements} The authors are grateful to two anonymous referees for their careful reading and their valuable comments. We thank Michele Benzi and Andy Wathen for their helpful feedback. JWP gratefully acknowledges support from the Engineering and Physical Sciences Research Council (UK) grant EP/S027785/1.

\appendix

\section{Proof of Theorem \ref{thm:k=3}}\label{app:proof} 
We prove Theorem \ref{thm:k=3} by contradiction, through assuming that the eigenvalues $\mu$ are not contained in the intervals stated in the result. We define the notation $V \succ W$ to mean that $V-W$ is (symmetric) positive definite; similarly we define $V \succeq W$, $V \prec W$, and $V \preceq W$ to mean that $V-W$ is positive semi-definite, negative definite, and negative semi-definite, respectively.

\begin{figure}
\centering
\vspace{-9.5em}
\includegraphics[width=0.9\textwidth]{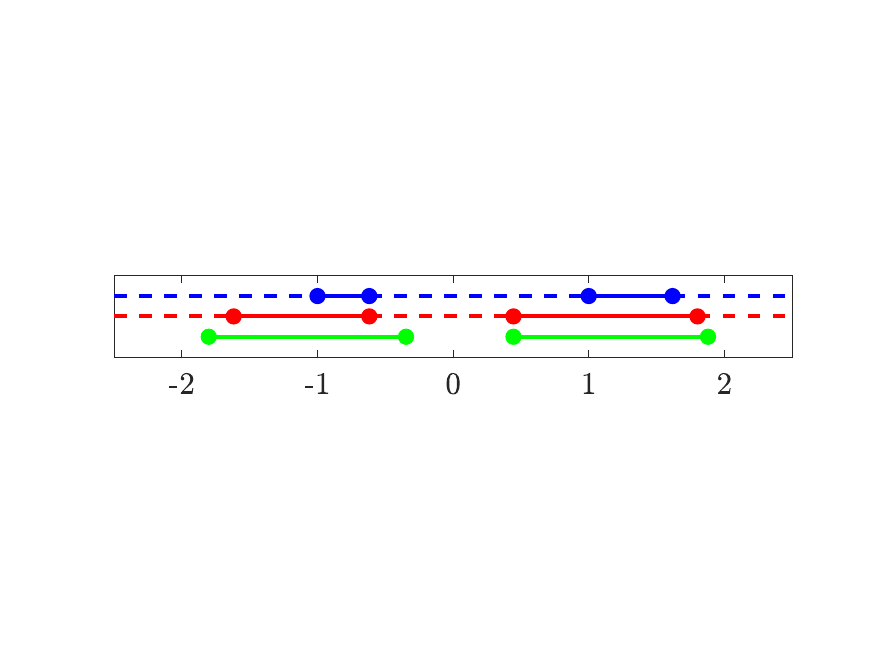}
\vspace{-9.5em}
\caption{Illustration of intervals within which we allow the possibilities of eigenvalues of $\mathcal{P}_{D,3}^{-1}\mathcal{A}_3$ belonging (bold lines), as we progress through the proof, as well as intervals still to be analysed (dotted lines). Intervals are given after analysing the second block row (top, blue), third block row (middle, red), and our conclusions after analysing the fourth block row (bottom, green).}
\label{appendixtheorem_plot}
\end{figure}

Following the structure of the proof in \cite[Theorem 5.2]{PearsonPotschka}, we examine the associated eigenproblem:
\begin{equation}
\ \label{4by4} \left(\begin{array}{cccc}
A_0 & B_1^T & O & O \\ B_1 & -A_1 & B_2^T & O \\ O & B_2 & A_2 & B_3^T \\ O & O & B_3 & -A_3 \\
\end{array}\right)\left(\begin{array}{c}
\mathbf{x}_1 \\ \mathbf{x}_2 \\ \mathbf{x}_3 \\ \mathbf{x}_4 \\
\end{array}\right)=\mu\left(\begin{array}{cccc}
A_0 & O & O & O \\ O & S_1 & O & O \\ O & O & S_2 & O \\ O & O & O & S_3 \\
\end{array}\right)\left(\begin{array}{c}
\mathbf{x}_1 \\ \mathbf{x}_2 \\ \mathbf{x}_3 \\ \mathbf{x}_4 \\
\end{array}\right).
\end{equation}
In the subsequent analysis, we progress through the second, third, and fourth block rows of \eqref{4by4}, establishing conditions upon which invertibility of appropriate matrices holds, and hence upon which we may advance to the analysis of the next block row. When the relevant conditions are not satisfied, we exclude certain intervals from the analysis and allow the possibility that eigenvalues of $\mathcal{P}_{D,3}^{-1}\mathcal{A}_3$ are contained within these intervals. Finally, when analysing the fourth block row, we exclude the possibility of eigenvalues belonging to particular intervals by contradiction. See Figure \ref{appendixtheorem_plot} for an illustration of the deductions made as we progress through the block rows of \eqref{4by4}.

\underline{\textbf{(First and) Second block row:}} \, The first block row of \eqref{4by4} tells us that $p\mathbf{x}_1=-A_0^{-1} B_1^T \mathbf{x}_2$, where $p=1-\mu\neq0$ by assumption. Substituting this into the second block row of \eqref{4by4}, simple algebra gives that $\mathbf{x}_2=p[qB_1 A_0^{-1} B_1^T+rA_1]^{-1} B_2^T \mathbf{x}_3$ where $q=1+\mu-\mu^2$, $r=1-\mu^2$, if $\mu\in(-\infty,-1)\cup(\frac{1}{2}(1-\sqrt{5}),1)\cup(\frac{1}{2}(1+\sqrt{5}),+\infty)$ in which case $qB_1 A_0^{-1} B_1^T+rA_1$ is guaranteed to be invertible. We allow the possibility that eigenvalues are contained in $[-1,\frac{1}{2}(1-\sqrt{5})]\cup[1,\frac{1}{2}(1+\sqrt{5})]$, and hence exclude these intervals from the remaining analysis (see top line of Figure \ref{appendixtheorem_plot}).

\underline{\textbf{Third block row:}} \, When $qB_1 A_0^{-1} B_1^T+rA_1$ is invertible, the third block row of \eqref{4by4} tells us that
\begin{align*}
[-pA_2 +\Phi]\mathbf{x}_3 ={}&B_3^T \mathbf{x}_4, \\
\text{where}~~\Phi:={}&B_2 \left[ \frac{1}{\mu}B_1 A_0^{-1} B_1^T +\frac{1}{\mu}A_1 \right]^{-1} B_2^T -B_2 \left[ \frac{q}{p}B_1 A_0^{-1} B_1^T +\frac{r}{p}A_1 \right]^{-1} B_2^T,
\end{align*}
excepting the case $\mu=0$, as the decomposition $\mathcal{A}_k=\mathcal{P}_{L,k}\bar{\mathcal{P}}_{D,k}^{-1}\mathcal{P}_{U,k}$ guarantees the invertibility of $\mathcal{A}_k$. We now consider values of $\mu$ for which $-pA_2 +\Phi$ is invertible:
\begin{itemize}
\item If $\mu>2\,\text{cos}(\frac{\pi}{7})$, $-pA_2$ is positive semi-definite. As $p,q,r<0$, $\frac{1}{\mu}<\frac{q}{p}$, and $\frac{1}{\mu}<\frac{r}{p}$, we see that $\frac{1}{\mu}B_1 A_0^{-1} B_1^T +\frac{1}{\mu}A_1 \prec \frac{q}{p}B_1 A_0^{-1} B_1^T +\frac{r}{p}A_1$, hence $\Phi$ is positive semi-definite. Hence, $-pA_2 +\Phi$ is positive definite: any vector belonging to the nullspaces of $-pA_2$ and $\Phi$, contradicts the positive definiteness of $S_2$.
\item If $\mu\in(0,2\,\text{cos}(\frac{3\pi}{7}))$, then $p,q,r>0$, $\frac{1}{\mu}>\frac{q}{p}$, and $\frac{1}{\mu}>\frac{r}{p}$. It follows that $\Phi$ is negative semi-definite, as is $-pA_2$. Hence, $-pA_2 +\Phi$ is negative definite.
\item If $\mu\in(\frac{1}{2}(1-\sqrt{5}),0)$, $-pA_2$ is negative semi-definite. As $\frac{1}{\mu}<0$, $\frac{q}{p},\frac{r}{p}>0$, it follows that $\Phi$ is negative definite, hence so is $-pA_2 +\Phi$.
\item If $\mu<-\frac{1}{2}(1+\sqrt{5})$, then $p>0$, $q,r<0$, $\frac{1}{\mu}>\frac{q}{p}$, and $\frac{1}{\mu}>\frac{r}{p}$. It follows that $\Phi$ is negative semi-definite, as is $-pA_2$. Hence $-pA_2 +\Phi$ is negative definite.
\end{itemize}
Hence $\mathbf{x}_3 =[-pA_2 +\Phi]^{-1} B_3^T \mathbf{x}_4$ when $\mu\in(-\infty,-\frac{1}{2}(1+\sqrt{5}))\cup(\frac{1}{2}(1-\sqrt{5}),0)\cup(0,2\,\text{cos}(\frac{3\pi}{7}))\cup(\text{cos}(\frac{\pi}{7}),+\infty)$. We now exclude the intervals $[-\frac{1}{2}(1+\sqrt{5}),\frac{1}{2}(1-\sqrt{5})]\cup[2\,\text{cos}(\frac{3\pi}{7}),2\,\text{cos}(\frac{\pi}{7})]$ from the remaining analysis, accepting that eigenvalues may be contained within the intervals (see middle line of Figure \ref{appendixtheorem_plot}).

\underline{\textbf{Fourth block row:}} \, Within the remaining intervals under consideration, we may substitute $\mathbf{x}_3 =[-pA_2 +\Phi]^{-1} B_3^T \mathbf{x}_4$ into the fourth block row of \eqref{4by4} to give:
\begin{align}
&(\mu+1)\mathbf{x}_4^T A_3 \mathbf{x}_4 \label{4thBlock} \\
&\quad=\mathbf{x}_4^T \left(B_3 [-pA_2+\Phi]^{-1} B_3^T -B_3 \left[ \frac{1}{\mu}A_2 +B_2 [\mu B_1 A_0^{-1} B_1^T +\mu A_1]^{-1} B_2^T \right]^{-1} B_3^T\right)\mathbf{x}_4. \nonumber
\end{align}
The left-hand side of \eqref{4thBlock} is non-negative if $\mu>-1$ and non-positive if $\mu<-1$; we use this to contradict the presence of eigenvalues $\mu$ in certain intervals (note that the possibility $\mu=0$ is already excluded):
\begin{itemize}
\item \underline{$\mu>2\,\text{cos}(\frac{\pi}{9})$:} To exclude this interval we may show that the right-hand side of \eqref{4thBlock} is non-positive: if it is zero we may argue that $\mathbf{x}_4 \in \text{null}(B_3^T)$ and the problem reduces to the case $k=2$ (see \cite[Theorem 4]{BrGr21} and \cite[Theorem 5.3]{PearsonPotschka}), and if it is negative we have a contradiction. In this regime we have that $-pA_2 \succeq \frac{1}{\mu}A_2$; further, as $p<0$, $r<q<0$, and $\frac{1}{\mu}+\frac{p}{q}<\mu$, we have that
\begin{align*}
&B_2 \left[ \mu B_1 A_0^{-1} B_1^T +\mu A_1 \right]^{-1} B_2^T +B_2 \left[ \frac{q}{p}B_1 A_0^{-1} B_1^T +\frac{r}{p}A_1 \right]^{-1} B_2^T \\
&\preceq \left(\frac{1}{\mu}+\frac{p}{q}\right)B_2 \left[ B_1 A_0^{-1} B_1^T +A_1 \right]^{-1} B_2^T \preceq B_2 \left[ \frac{1}{\mu}B_1 A_0^{-1} B_1^T +\frac{1}{\mu}A_1 \right]^{-1} B_2^T,
\end{align*}
and hence $\Phi \succeq B_2 [\mu B_1 A_0^{-1} B_1^T +\mu A_1]^{-1} B_2^T$. These facts, and the result that $W \prec V \Leftrightarrow V^{-1} \prec W^{-1}$, for symmetric positive definite $V$, $W$ \cite[Theorem 4.5]{FiedlerPtak97}, tell us that the right-hand side of \eqref{4thBlock} is non-positive, as required.
\item \underline{$\mu\in(0,2\,\text{cos}(\frac{3\pi}{7}))$:} Here we have that $-pA_2 \preceq0$, $\frac{1}{\mu}A_2 \succeq0$, $\Phi\preceq0$, and $B_2 \left[ \mu B_1 A_0^{-1} B_1^T +\mu A_1 \right]^{-1} B_2^T \succeq0$. We thus see the right-hand side of \eqref{4thBlock} is non-positive, so exclude this interval for $\mu$ by the same logic as above.
\item \underline{$\mu\in(2\,\text{cos}(\frac{5\pi}{9}),0)$:} We have that $-pA_2 \succeq \frac{1}{\mu}A_2$. Further, as $p>0$, $r>q>0$, and $\mu-\frac{p}{q}>\frac{1}{\mu}$, we also write
\begin{align*}
&-\left(B_2 \left[ -\frac{1}{\mu} B_1 A_0^{-1} B_1^T -\frac{1}{\mu} A_1 \right]^{-1} B_2^T +B_2 \left[ \frac{q}{p}B_1 A_0^{-1} B_1^T +\frac{r}{p}A_1 \right]^{-1} B_2^T\right) \\
&\succeq \left(\mu-\frac{p}{q}\right)B_2 \left[ B_1 A_0^{-1} B_1^T +A_1 \right]^{-1} B_2^T \succeq B_2 \left[ \mu B_1 A_0^{-1} B_1^T +\mu A_1 \right]^{-1} B_2^T.
\end{align*}
Hence $\Phi \succeq B_2 [\mu B_1 A_0^{-1} B_1^T +\mu A_1]^{-1} B_2^T$, and so the right-hand side of \eqref{4thBlock} is non-positive, allowing us to exclude this interval by the logic above.
\item \underline{$\mu<-2\,\text{cos}(\frac{\pi}{7})$:} We have that $-pA_2 \preceq \frac{1}{\mu}A_2$. Further, as $p>0$, $q<r<0$, and $\frac{1}{\mu}+\frac{p}{r}>\mu$, we also write
\begin{align*}
&B_2 \left[ \mu B_1 A_0^{-1} B_1^T +\mu A_1 \right]^{-1} B_2^T +B_2 \left[ \frac{q}{p}B_1 A_0^{-1} B_1^T +\frac{r}{p}A_1 \right]^{-1} B_2^T \\
&\succeq \left(\frac{1}{\mu}+\frac{p}{r}\right)B_2 \left[ B_1 A_0^{-1} B_1^T +A_1 \right]^{-1} B_2^T \succeq B_2 \left[ \frac{1}{\mu}B_1 A_0^{-1} B_1^T +\frac{1}{\mu}A_1 \right]^{-1} B_2^T.
\end{align*}
Hence $\Phi \preceq B_2 [\mu B_1 A_0^{-1} B_1^T +\mu A_1]^{-1} B_2^T$, and so the right-hand side of \eqref{4thBlock} is non-negative, allowing us to exclude this interval by the logic above.
\end{itemize}
The result is thus proved by contradiction, through excluding the remaining intervals.


\begin{thebibliography}{99}

\bibitem{AliBeikBenzi} F. P. Ali Beik and M. Benzi, \emph{Iterative methods for double saddle point systems}, SIAM J. Matrix Anal. Appl. \textbf{39}, 902--921 (2018)

\bibitem{AxNe06} O. Axelsson and M. Neytcheva, \emph{Eigenvalue estimates for preconditioned saddle point matrices}, Numer. Linear Alg. Appl. \textbf{13}, 339--360 (2006)

\bibitem{BSZ} A. Beigl, J. Sogn, and W. Zulehner, \emph{Robust preconditioners for multiple saddle point problems and applications to optimal control problems}, SIAM J. Matrix Anal. Appl. \textbf{41}, 1590--1615 (2020)

\bibitem{HSL_MI20} J. Boyle, M. Mihajlovi\'{c}, and J. Scott, \emph{\texttt{HSL\_MI20}: an efficient AMG preconditioner for finite element problems in 3D}, Int. J. Numer. Meth. Eng. \textbf{82}, 64--98 (2010)

\bibitem{HSL_MI20_code} \emph{\texttt{HSL\_MI20} unsymmetric system: algebraic multigrid preconditioner}, code available at \url{https://www.hsl.rl.ac.uk/catalogue/hsl_mi20.html}

\bibitem{BrGr21} S. Bradley and C. Greif, \emph{Eigenvalue bounds for double saddle-point systems}, arXiv preprint \texttt{arXiv:2110.13328} (2021)

\bibitem{CJL22} M. Cai, G. Ju, and J. Li, \emph{Schur complement based preconditioners for twofold and block tridiagonal saddle point problems}, arXiv preprint \texttt{arXiv:2108.08332} (2022)

\bibitem{FiedlerPtak97} M. Fiedler and V. Pt\'{a}k, \emph{A new positive definite geometric mean of two positive definite matrices}, Linear Algebra Appl. \textbf{251}, 1--20 (1997)

\bibitem{FoSa12} D. C.-L. Fong and M. Saunders, \emph{CG versus MINRES: An empirical comparison}, SJU Journal for Science \textbf{17}, 44--62 (2012)

\bibitem{GMPS} P. E. Gill, W. Murray, D. B. Poncele\'{o}n, and M. A. Saunders, \emph{Preconditioners for indefinite systems arising in optimization}, SIAM J. Matrix Anal. Appl. \textbf{13}, 292--311 (1992)

\bibitem{GVI} G. H. Golub and R. S. Varga, \emph{Chebyshev semi-iterative methods, successive over-relaxation iterative methods, and second order Richardson iterative methods, Part I}, Numer. Math. \textbf{3}, 147--156 (1961)

\bibitem{GVII} G. H. Golub and R. S. Varga, \emph{Chebyshev semi-iterative methods, successive over-relaxation iterative methods, and second order Richardson iterative methods, Part II}, Numer. Math. \textbf{3}, 157--168 (1961)

\bibitem{skfem2020} T. Gustafsson and G. D. McBain, \emph{\texttt{scikit-fem}: a Python package for finite element assembly}, J. Open Source Softw. \textbf{5}, 2369 (2020)

\bibitem{ItKu2003} K. Ito and K. Kunisch, \emph{Semi-smooth Newton methods for state-constrained optimal control problems}, Systems Control Lett. \textbf{50}, {221--228} (2003)

\bibitem{Kuznetsov} Y. A. Kuznetsov, \emph{Efficient iterative solvers for elliptic finite element problems on nonmatching grids}, Russ. J. Numer. Anal. Math. M. \textbf{10}, 187--211 (1995)

\bibitem{MNN} K.-A. Mardal, B. F. Nielsen, and M. Nordaas, \emph{Robust preconditioners for PDE-constrained optimization with limited observations}, BIT Numer. Math. \textbf{57}, 405--431 (2017)

\bibitem{MGW} M. F. Murphy, G. H. Golub, and A. J. Wathen, \emph{A note on preconditioning for indefinite linear systems}, SIAM J. Sci. Comput. \textbf{21}, 1969--1972 (2000)





\bibitem{minres} C. C. Paige and M. A. Saunders, \emph{Solution of sparse indefinite systems of linear equations}, SIAM J. Numer. Anal. \textbf{12}, 617--629 (1975)

\bibitem{minres_code} C. C. Paige, M. A. Saunders, et al., \emph{MINRES: sparse symmetric equations}, code available at \url{https://web.stanford.edu/group/SOL/software/minres/}

\bibitem{PW12} J. W. Pearson and A. J. Wathen, \emph{A new approximation of the Schur complement in preconditioners for PDE-constrained optimization}, Numer. Linear Alg. Appl. \textbf{19}, 816--829 (2012)

\bibitem{PW13} J. W. Pearson and A. J. Wathen, \emph{Fast iterative solvers for convection--diffusion control problems}, Electron. Trans. Numer. Anal. \textbf{40}, 294--310 (2013)

\bibitem{PearsonDPhil} J. Pearson, \emph{Fast Iterative Solvers for PDE-Constrained Optimization Problems}, DPhil thesis, University of Oxford (2013)

\bibitem{PearsonPotschka} J. W. Pearson and A. Potschka, \emph{A preconditioned inexact active-set method for large-scale nonlinear optimal control problems}, arXiv preprint \texttt{arXiv:2112.05020} (2021)

\bibitem{PoBo21} A. Potschka and H. G. Bock, \emph{A sequential homotopy method for mathematical programming problems}, Math. Program., \textbf{187}, 459--486 (2021)

\bibitem{gmres} Y. Saad and M. H. Schultz, \emph{GMRES: a generalized minimal residual algorithm for solving nonsymmetric linear systems}, SIAM J. Sci. Stat. Comput. \textbf{7}, 856--869 (1986)

\bibitem{SiWa94} D. Silvester and A. Wathen, \emph{Fast iterative solution of stabilised Stokes systems. Part II: using general block preconditioners}, SIAM J. Numer. Anal. \textbf{31}, 1352--1367 (1994)

\bibitem{SognZulehner} J. Sogn and W. Zulehner, \emph{Schur complement preconditioners for multiple saddle point problems of block tridiagonal form with application to optimization problems}, IMA J. Numer. Anal. \textbf{39}, 1328--1359 (2019)

\bibitem{Ulbrich2011} M. Ulbrich, \emph{Semismooth Newton Methods for Variational Inequalities and Constrained Optimization Problems in Function Spaces}, MOS-SIAM Series on Optimization, Vol. 11., SIAM, Philadelphia, PA (2011)

\bibitem{WathenRees} A. Wathen and T. Rees, \emph{Chebyshev semi-iteration in preconditioning for problems including the mass matrix}, Electron. Trans. Numer. Anal. \textbf{34}, 125--135 (2009)

\bibitem{WathenEigs} A. J. Wathen, \emph{Realistic eigenvalue bounds for the Galerkin mass matrix}, IMA J. Numer. Anal. \textbf{7}, 449--457 (1987)

\end{thebibliography}
\end{document}